\documentclass[a4paper,11pt]{amsart}
\usepackage{amssymb}
\usepackage[arrow,matrix]{xy}

\title[Canonical algebras derived equivalent to incidence algebras]
{Which canonical algebras are derived equivalent to incidence algebras
of posets?}

\author{Sefi Ladkani}

\address{Einstein Institute of Mathematics, The Hebrew University of Jerusalem, Jerusalem 91904, Israel}
\email{sefil@math.huji.ac.il}

\DeclareMathOperator{\cone}{C}
\DeclareMathOperator{\Ext}{Ext}
\DeclareMathOperator{\Hom}{Hom}
\DeclareMathOperator{\HH}{HH}
\DeclareMathOperator{\gldim}{gl.\!dim}

\newcommand{\bp}{\mathbf{p}}
\newcommand{\bl}{\pmb{\lambda}}

\newcommand{\gL}{\Lambda}
\newcommand{\cC}{\mathcal{C}}
\newcommand{\cD}{\mathcal{D}}
\newcommand{\cH}{\mathcal{H}}
\newcommand{\cF}{\mathcal{F}}
\newcommand{\bZ}{\mathbb{Z}}

\newcommand{\cL}{\cC^b(\gL)}
\newcommand{\dL}{\cD^b(\gL)}

\newcommand{\dX}{\cD^b(k X)}

\theoremstyle{plain}
\newtheorem{theorem}{Theorem}[section]
\newtheorem*{theorem*}{Theorem}
\newtheorem{lemma}[theorem]{Lemma}
\newtheorem{prop}[theorem]{Proposition}
\newtheorem{cor}[theorem]{Corollary}

\theoremstyle{definition}

\newtheorem{example}[theorem]{Example}
\newtheorem{rem}[theorem]{Remark}

\numberwithin{equation}{section}

\begin{document}
\begin{abstract}
We give a full description of all the canonical algebras over an
algebraically closed field that are derived equivalent to incidence
algebras of finite posets. These are the canonical algebras whose
number of weights is either 2 or 3.
\end{abstract}

\maketitle

This note concerns the characterization of the canonical algebras over
an algebraically closed field that are derived equivalent to incidence
algebras of finite partially ordered sets (posets), expressed in the
following theorem.

\begin{theorem*}
Let $\gL$ be a canonical algebra of type $(\bp,\bl)$ over an
algebraically closed field. Then $\gL$ is derived equivalent to an
incidence algebra of a poset if and only if the number of weights of
$\bp$ is either 2 or 3.
\end{theorem*}

This theorem can be interpreted both geometrically and algebraically.
From a geometric viewpoint, by considering modules over incidence
algebras as sheaves over finite spaces~\cite{Ladkani07} and using the
derived equivalence between the categories of modules over a canonical
algebra and coherent sheaves over a weighted projective
line~\cite{GeigleLenzing87}, we are able to obtain explicit derived
equivalences between the categories of sheaves of finite dimensional
vector spaces over certain finite $T_0$ topological spaces and the
categories of coherent sheaves over certain weighted projective lines.

From an algebraic viewpoint, in an attempt to classify all piecewise
hereditary incidence algebras over an algebraically closed field, one
first asks which types of piecewise hereditary categories can actually
occur. Happel's classification~\cite{Happel01} tells us that we only
need to consider the canonical algebras and path algebras of quivers.
For the canonical algebras the theorem above gives a complete answer,
while for path algebras, see the remarks in Section~\ref{ssec:quivers}.

We finally note that for the constructions of incidence algebras
derived equivalent to canonical algebras, the assumption that the base
field is algebraically closed can be omitted.

\subsection*{Acknowledgement}
I would like to thank Helmut Lenzing for useful discussions related to
this paper, and the referee for the helpful comments.

\setcounter{section}{-1}
\section{Notations}
%%%%%%%%%%%%%%%%%%%

\sloppy The canonical algebras were introduced in~\cite{Ringel84}. Let
$k$ be a field, $\bp=(p_1,\dots,p_t)$ be a sequence of $t \geq 2$
positive integers (\emph{weights}), and $\bl=(\lambda_3, \dots,
\lambda_t)$ be a sequence of pairwise distinct elements of $k \setminus
\{0\}$. The \emph{canonical algebra of type $(\bp,\bl)$}, is the
algebra $\gL(\bp,\bl) = kQ/I$ where $Q$ is the quiver
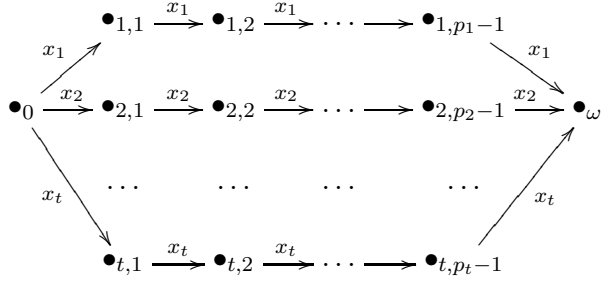
\[
\xymatrix@=1.5pc{
& {\bullet_{1,1}} \ar^{x_1}[r] & {\bullet_{1,2}} \ar^{x_1}[r]
& {\ldots} \ar[r] & {\bullet_{1,p_1-1}} \ar^{x_1}[rd] \\
{\bullet_0} \ar^{x_1}[ur] \ar_{x_t}[rdd] \ar^{x_2}[r]
& {\bullet_{2,1}} \ar^{x_2}[r] & {\bullet_{2,2}} \ar^{x_2}[r]
& {\ldots} \ar[r] & {\bullet_{2,p_2-1}} \ar^{x_2}[r]
& {\bullet_{\omega}} \\
& {\ldots} & {\ldots} & {\ldots} & {\ldots} \\
& {\bullet_{t,1}} \ar^{x_t}[r] & {\bullet_{t,2}} \ar^{x_t}[r]
& {\ldots} \ar[r] & {\bullet_{t,p_t-1}} \ar_{x_t}[ruu] }
\]
and $I$ is the ideal in the path algebra $kQ$ generated by the
following linear combinations of paths from $0$ to $\omega$:
\[
I = \left\langle x_i^{p_i} - x_2^{p_2} + \lambda_i x_1^{p_1} \,:\, 3
\leq i \leq t \right\rangle
\]

As noted in~\cite{Ringel84}, as long as $t \geq 3$, one can omit
weights equal to 1, and when $t = 2$, the ideal $I$ vanishes and the
canonical algebra is equal to the path algebra of $Q$. In the latter
case, one usually writes only the weights greater than 1, in particular
the algebra of type $()$ equals the path algebra of the Kronecker
quiver. Hence when speaking on the \emph{number} of weights, we shall
always mean the number of $p_i$ with $p_i \geq 2$.

Let $X$ be a finite partially ordered set (\emph{poset}). The
\emph{incidence algebra} $kX$ is the $k$-algebra spanned by the
elements $e_{xy}$ for the pairs $x \leq y$ in $X$, with the
multiplication defined by setting $e_{xy} e_{y'z} = e_{xz}$ if $y=y'$
and $e_{xy} e_{y'z} = 0$ otherwise.

A \emph{$k$-diagram} $\cF$ over $X$ consists of finite dimensional
vector spaces $\cF(x)$ for $x \in X$, together with linear
transformations $r_{xx'} : \cF(x) \to \cF(x')$ for all $x \leq x'$,
satisfying the conditions $r_{xx} = 1_{\cF(x)}$ and $r_{xx''} =
r_{x'x''}r_{xx'}$ for all $x \leq x' \leq x''$. The category of finite
dimensional right modules over $kX$ is equivalent to the category of
$k$-diagrams over $X$, see~\cite{Ladkani07}.

For a finite-dimensional algebra $\gL$ over $k$, we denote by $\cL$ the
category of bounded complexes of (right) finite-dimensional
$\gL$-modules, and by $\dL$ the bounded derived category. The algebra
$\gL$ is \emph{piecewise hereditary} if $\dL$ is equivalent as
triangulated category to $\cD^b(\cH)$ for a hereditary $k$-category
$\cH$.

\section{The necessity of the condition $t \leq 3$ in the Theorem}
%%%%%%%%%%%%%%%%%%%%%%%%%%%%%%%%%%%%%%%%%%%%%%%%%%%%%%%%%%%%%%%%%%

For a $k$-algebra $\gL$, denote by $\HH^i(\gL)$ the \emph{$i$-th
Hochschild cohomology} of $\gL$, which equals $\Ext^i_{\gL \otimes
\gL^{op}}(\gL, \gL)$, where $\gL$ is considered as a
$\gL$-$\gL$-bimodule in the natural way.

\begin{prop}
\label{p:HHpiece} Let $X$ be a poset such that $kX$ is piecewise
hereditary. Then $\HH^i(kX) = 0$ for any $i>1$.
\end{prop}
\begin{proof}
First, since for any two finite-dimensional $k$-algebras $\gL_1$,
$\gL_2$ we have $\HH^i(\gL_1 \oplus \gL_2) = \HH^i(\gL_1) \oplus
\HH^i(\gL_2)$, we may assume that $X$ is connected, as the
decomposition of $X$ into connected components $X = \sqcup_{i=1}^{r}
X_i$ induces a decomposition of the incidence algebra $kX =
\bigoplus_{i=1}^{r} kX_i$.

Let $k_X$ be the \emph{constant diagram} on $X$, defined by $k_X(x) =
k$ for all $x \in X$ with all maps being the identity on $k$.
By~\cite[Corollary~3.20]{Ladkani07}, $\HH^i(kX) = \Ext^i_X(k_X, k_X)$.

Since $X$ is connected, $k_X$ is indecomposable, and by~\cite[IV
(1.9)]{Happel88}, the groups $\Ext_X^i(k_X, k_X)$ vanish for $i > 1$,
hence $\HH^i(kX) = 0$ for $i > 1$.
\end{proof}

\begin{cor}
\label{c:canon} Let $k$ be algebraically closed and let $\gL$ be a
canonical algebra over $k$ of type $(\bp,\bl)$ where $\bp =
(p_1,\dots,p_t)$. If $\dL \simeq \dX$ for some poset $X$, then $t \leq
3$.
\end{cor}
\begin{proof}
Assume that $t \geq 4$. Then $\gL$ is not of domestic type and
by~\cite[Theorem~2.4]{Happel98}, $\dim_k \HH^2(kX) = t-3$, a
contradiction to Proposition~\ref{p:HHpiece}. Therefore $t \leq 3$.
\end{proof}

\section{Constructions of posets from canonical algebras}
%%%%%%%%%%%%%%%%%%%%%%%%%%%%%%%%%%%%%%%%%%%%%%%%%%%%%%%%%

\subsection{The case $t=3$}
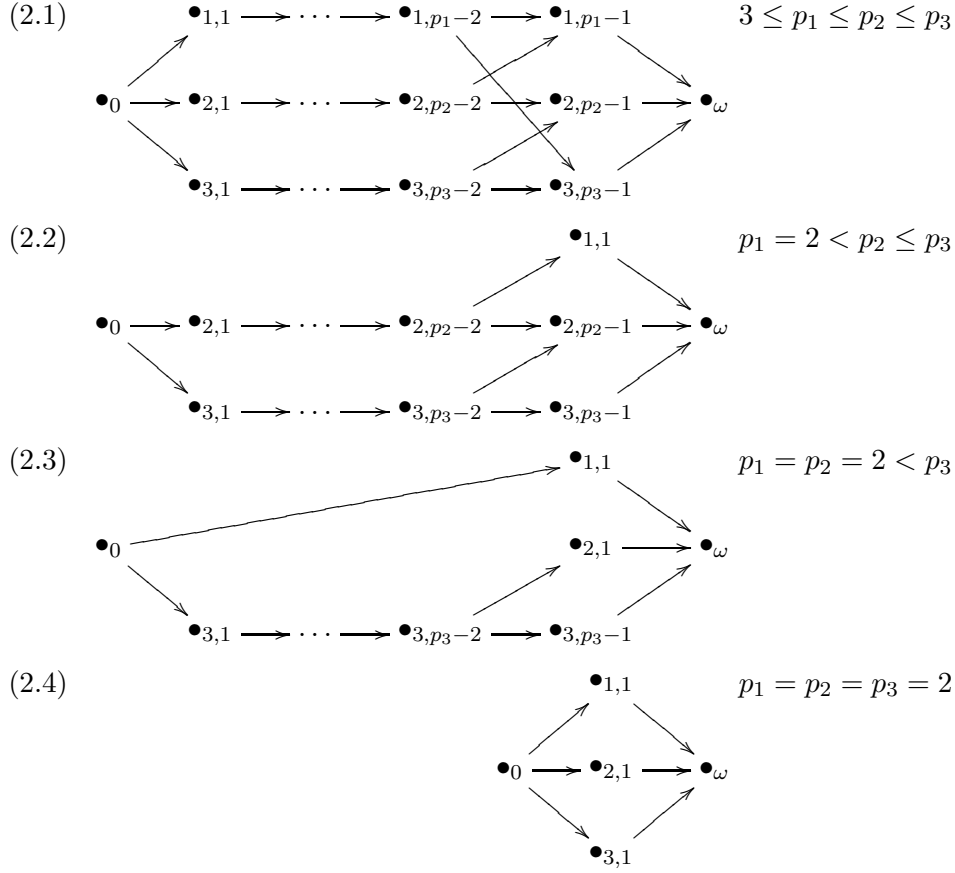
\begin{figure}
\begin{align} \label{e:333}
\xymatrix@=1.5pc{
& {\bullet_{1,1}} \ar[r] & {\ldots} \ar[r] &
{\bullet_{1,p_1-2}} \ar[r] \ar[ddr] & {\bullet_{1,p_1-1}} \ar[dr] \\
{\bullet_0} \ar[ur] \ar[r] \ar[dr] &
{\bullet_{2,1}} \ar[r] & {\ldots} \ar[r] &
{\bullet_{2,p_2-2}} \ar[r] \ar[ur] & {\bullet_{2,p_2-1}} \ar[r] &
{\bullet_{\omega}} \\
& {\bullet_{3,1}} \ar[r] & {\ldots} \ar[r] &
{\bullet_{3,p_3-2}} \ar[r] \ar[ur] & {\bullet_{3,p_3-1}} \ar[ur]
}
& 3 \leq p_1 \leq p_2 \leq p_3
\\ \label{e:233}
\xymatrix@=1.5pc{
& & & & {\bullet_{1,1}} \ar[dr] \\
{\bullet_0} \ar[r] \ar[dr] &
{\bullet_{2,1}} \ar[r] & {\ldots} \ar[r] &
{\bullet_{2,p_2-2}} \ar[r] \ar[ur] & {\bullet_{2,p_2-1}} \ar[r] &
{\bullet_{\omega}} \\
& {\bullet_{3,1}} \ar[r] & {\ldots} \ar[r] &
{\bullet_{3,p_3-2}} \ar[r] \ar[ur] & {\bullet_{3,p_3-1}} \ar[ur]
}
& p_1 = 2 < p_2 \leq p_3
\\ \label{e:223}
\xymatrix@=1.5pc{
& & & & {\bullet_{1,1}} \ar[dr] \\
{\bullet_0} \ar[urrrr] \ar[dr] & & & & {\bullet_{2,1}} \ar[r] &
{\bullet_{\omega}} \\
& {\bullet_{3,1}} \ar[r] & {\ldots} \ar[r] &
{\bullet_{3,p_3-2}} \ar[r] \ar[ur] & {\bullet_{3,p_3-1}} \ar[ur]
}
& p_1 = p_2 = 2 < p_3
\\ \label{e:222}
\xymatrix@=1.5pc{
&  {\bullet_{1,1}} \ar[dr] \\
{\bullet_0} \ar[ur] \ar[r] \ar[dr] & {\bullet_{2,1}} \ar[r] &
{\bullet_{\omega}} \\
& {\bullet_{3,1}} \ar[ur]
}
& p_1 = p_2 = p_3 = 2
\end{align}
\caption{Hasse diagrams of posets derived equivalent to canonical
algebras of type $(p_1, p_2, p_3)$.} \label{f:posets}
\end{figure}

%%%%%%%%%%%%%%%%%%%%%%%%%%%
Recall that when $t=3$, the canonical algebra $\gL(\bp,\bl)$ is
independent of the parameter $\lambda_3$, so we may assume that
$\lambda_3 = 1$, and denote the algebra by $\gL(\bp)$. Let $\bp = (p_1,
p_2, p_3)$ be a triplet of weights with $2 \leq p_1 \leq p_2 \leq p_3$.
Attach to $\bp$ a poset $X_{\bp}$ whose Hasse diagram is shown in
Figure~\ref{f:posets}. Explicitly, use the Hasse diagram
of~\eqref{e:333} if $p_1 > 2$, \eqref{e:233} if $p_1 = 2$ and $p_2 >
2$, \eqref{e:223} if $p_2 = 2$ and $p_3 > 2$ and~\eqref{e:222} if $p_3
= 2$.

\begin{theorem} \label{t:Xp}
Let $k$ be a field and let $\bp = (p_1, p_2, p_3)$. Then $\gL(\bp)$ is
derived equivalent to the incidence algebra of $X_{\bp}$.
\end{theorem}

\begin{proof}
The idea of the proof relies on the notion of a formula introduced
in~\cite{Ladkani07b}, and we refer to that paper for more details. We
shall first construct a functor $F : \cC^b(k X_{\bp}) \to
\cC^b(\gL(\bp))$ that induces a triangulated functor $\widetilde{F} :
\cD^b(k X_{\bp}) \to \cD^b(\gL(\bp))$, and then prove that
$\widetilde{F}$ is an equivalence.

Note that similarly to the identification in~\cite{Ladkani07b} of
complexes of diagrams with diagrams of complexes, we may identify
complexes of modules over the canonical algebra with a
(non-commutative) diagram of complexes of vector spaces satisfying the
canonical algebra relations.

For a morphism $f : K \to L$ of complexes $K=(K^i, d_K^i)$, $L=(L^i,
d_L^i)$ of vector spaces, denote by $\cone(K \xrightarrow{f} L)$ the
\emph{cone} of $f$. Recall that $\cone(K \to L)^i = K^{i+1} \oplus
L^i$, with the differential acting as the matrix
\[
\begin{pmatrix}
-d_K^{i+1} & 0 \\ f^{i+1} & d_L^i
\end{pmatrix}
\]
by viewing the terms as column vectors of length two. Denote by $[1]$
the shift operator, that is, $K[1]^i=K^{i+1}$ with $d_{K[1]}^i =
-d_K^{i+1}$.

We will demonstrate the construction of $\widetilde{F}$ for the posets
of type~\eqref{e:333}, the other cases being similar. Let $K_0$,
$K_\omega$ and $K_{i,j}$ for $1 \leq i \leq 3$, $1 \leq j < p_i$ be
complexes in a commutative diagram
\[
\xymatrix@=1.5pc{
& {K_{1,1}} \ar[r] & {\ldots} \ar[r] &
{K_{1,p_1-2}} \ar^{x_{11}}[r] \ar^(0.8){x_{13}}[ddr] & {K_{1,p_1-1}} \ar^{y_1}[dr] \\
{K_0} \ar[ur] \ar[r] \ar[dr] &
{K_{2,1}} \ar[r] & {\ldots} \ar[r] &
{K_{2,p_2-2}} \ar_(0.4){x_{22}}[r] \ar^(0.6){x_{21}}[ur] & {K_{2,p_2-1}} \ar^{y_2}[r] &
{K_{\omega}} \\
& {K_{3,1}} \ar[r] & {\ldots} \ar[r] &
{K_{3,p_3-2}} \ar_{x_{33}}[r] \ar_(0.4){x_{32}}[ur] & {K_{3,p_3-1}} \ar_{y_3}[ur]
}
\]

Let $L_0 = K_0$ and $L_{i,j} = K_{i,j}$ for $1 \leq i \leq 3$ and $1
\leq j < p_i - 1$. Define
\begin{align*}
L_{1,p_1-1} &= \cone(K_{1,p_1-1} \oplus K_{3,p_3-1}
\xrightarrow{\left(\begin{smallmatrix} y_1 & y_3 \end{smallmatrix}\right)}
K_{\omega})[-1] \\
L_{2,p_2-1} &= \cone(K_{2,p_2-1} \oplus K_{1,p_1-1}
\xrightarrow{\left(\begin{smallmatrix} y_2 & y_1 \end{smallmatrix}\right)}
K_{\omega})[-1] \\
L_{3,p_3-1} &= \cone(K_{3,p_3-1} \oplus K_{2,p_2-1}
\xrightarrow{\left(\begin{smallmatrix} y_3 & y_2 \end{smallmatrix}\right)}
K_{\omega})[-1] \\
L_{\omega}  &= \cone(K_{1,p_1-1} \oplus K_{2,p_2-1} \oplus K_{3,p_3-1}
\xrightarrow{\left(\begin{smallmatrix} y_1 & y_2 & y_3 \end{smallmatrix}\right)}
K_{\omega})[-1]
\end{align*}
with the three maps
\begin{align*}
L_{1,p_1-2} \xrightarrow{\left(\begin{smallmatrix} x_{11} & -x_{13} & 0
\end{smallmatrix}\right)^{T}} L_{1,p_1-1} \\
L_{2,p_2-2} \xrightarrow{\left(\begin{smallmatrix} -x_{22} & x_{21} & 0
\end{smallmatrix}\right)^{T}} L_{2,p_2-1} \\
L_{3,p_3-2} \xrightarrow{\left(\begin{smallmatrix} x_{33} & -x_{32} & 0
\end{smallmatrix}\right)^{T}} L_{3,p_3-1}
\end{align*}
and the three maps from $L_{i,p_i-1}$ to $L_\omega$ being the canonical
embeddings.

Then the following is a (non-commutative) diagram of complexes
\[
\xymatrix@=1.5pc{
& {L_{1,1}} \ar[r] & {\ldots} \ar[r] &
{L_{1,p_1-2}} \ar[r] & {L_{1,p_1-1}} \ar[dr] \\
{L_0} \ar[ur] \ar[r] \ar[dr] &
{L_{2,1}} \ar[r] & {\ldots} \ar[r] &
{L_{2,p_2-2}} \ar[r] & {L_{2,p_2-1}} \ar[r] &
{L_{\omega}} \\
& {L_{3,1}} \ar[r] & {\ldots} \ar[r] &
{L_{3,p_3-2}} \ar[r] & {L_{3,p_3-1}} \ar[ur]
}
\]
that satisfies the canonical algebra relation, and we get the required
functor $F$ which induces, by the general considerations
in~\cite{Ladkani07b}, the functor $\widetilde{F}$.

To prove that $\widetilde{F}$ is an equivalence, we use Beilinson's
Lemma~\cite[Lemma~1]{Beilinson78} and verify that for any two simple
objects $S_x$, $S_y$ (where $x, y \in X$) and $i \in \bZ$, the functor
$\widetilde{F}$ satisfies
\[
\Hom_{\cD^b(k X_{\bp})}(S_x, S_y[i]) \simeq
\Hom_{\cD^b(\gL(\bp))}(\widetilde{F} S_x, \widetilde{F} S_y[i])
\]
and moreover the images $\widetilde{F} S_x$ generate $\cD^b(\gL(\bp))$
as a triangulated category.

We omit the details of this verification. However, we just mention that
$\widetilde{F} S_0$ and $\widetilde{F} S_{i,j}$, for $1 \leq i \leq 3$
and $1 \leq j < p_i - 1$, are the corresponding simple
$\gL(\bp)$-modules, while $\widetilde{F} S_{1,p_1-1}$, $\widetilde{F}
S_{2,p_2-1}$, $\widetilde{F} S_{3,p_3-1}$ and $\widetilde{F}
S_{\omega}$ are given by
\[
\begin{array}{cc}
\begin{array}{c}
\xymatrix@=1pc{
& {0} \ar[r] & {\ldots} \ar[r] &
{0} \ar[r] & {k} \ar[dr]^{1} \\
{0} \ar[ur] \ar[r] \ar[dr] &
{0} \ar[r] & {\ldots} \ar[r] &
{0} \ar[r] & {k} \ar[r]_{1} &
{k} \\
& {0} \ar[r] & {\ldots} \ar[r] &
{0} \ar[r] & {0} \ar[ur]
}
\end{array}
&
\begin{array}{c}
\xymatrix@=1pc{
& {0} \ar[r] & {\ldots} \ar[r] &
{0} \ar[r] & {0} \ar[dr] \\
{0} \ar[ur] \ar[r] \ar[dr] &
{0} \ar[r] & {\ldots} \ar[r] &
{0} \ar[r] & {k} \ar[r]_{1} &
{k} \\
& {0} \ar[r] & {\ldots} \ar[r] &
{0} \ar[r] & {k} \ar[ur]_{1}
}
\end{array}
\\
\widetilde{F} S_{1,p_1-1} &
\widetilde{F} S_{2,p_2-1}
\\
\\
\begin{array}{c}
\xymatrix@=1pc{
& {0} \ar[r] & {\ldots} \ar[r] &
{0} \ar[r] & {k} \ar[dr]^{1} \\
{0} \ar[ur] \ar[r] \ar[dr] &
{0} \ar[r] & {\ldots} \ar[r] &
{0} \ar[r] & {0} \ar[r] &
{k} \\
& {0} \ar[r] & {\ldots} \ar[r] &
{0} \ar[r] & {k} \ar[ur]_{1}
}
\end{array}
&
\Bigl(
\begin{array}{c}
\xymatrix@=1pc{
& {0} \ar[r] & {\ldots} \ar[r] &
{0} \ar[r] & {k} \ar[dr]^{1} \\
{0} \ar[ur] \ar[r] \ar[dr] &
{0} \ar[r] & {\ldots} \ar[r] &
{0} \ar[r] & {k} \ar[r]_{1} &
{k} \\
& {0} \ar[r] & {\ldots} \ar[r] &
{0} \ar[r] & {k} \ar[ur]_{1}
}
\end{array}
\Bigr)[-1]
\\
\widetilde{F} S_{3,p_3-1} &
\widetilde{F} S_\omega
\end{array}
\]
\end{proof}

\begin{example}[\protect{\cite[Example~18.6.2]{Lenzing99}}]
Let $A_2$ be the quiver $\xymatrix{{\bullet_1} \ar[r] & {\bullet_2}}$
and let $X = A_2 \times A_2 \times A_2$. Then the incidence algebra of
$X$ is derived equivalent to the canonical algebra of type $(3,3,3)$.
\end{example}

\begin{rem}
Observe that $\omega$ is the unique maximal element in the posets whose
Hasse diagrams are given in~\eqref{e:223} and~\eqref{e:222}. Hence by
taking $Y=\{\omega\}$ in Corollary~4.15 of~\cite{Ladkani07}, we recover
the fact that the canonical algebra of type $(2,2,p)$ is derived
equivalent to the path algebra of the extended Dynkin quiver
$\widetilde{D}_{p+2}$.
\end{rem}

\begin{rem}
Similar applications of~\cite[Theorem~1.1]{Ladkani07b} and its
corollaries for the posets in~\eqref{e:233} show that the canonical
algebra of type $(2,p_2,p_3)$ is derived equivalent to the incidence
algebras of the posets whose Hasse diagrams are given by
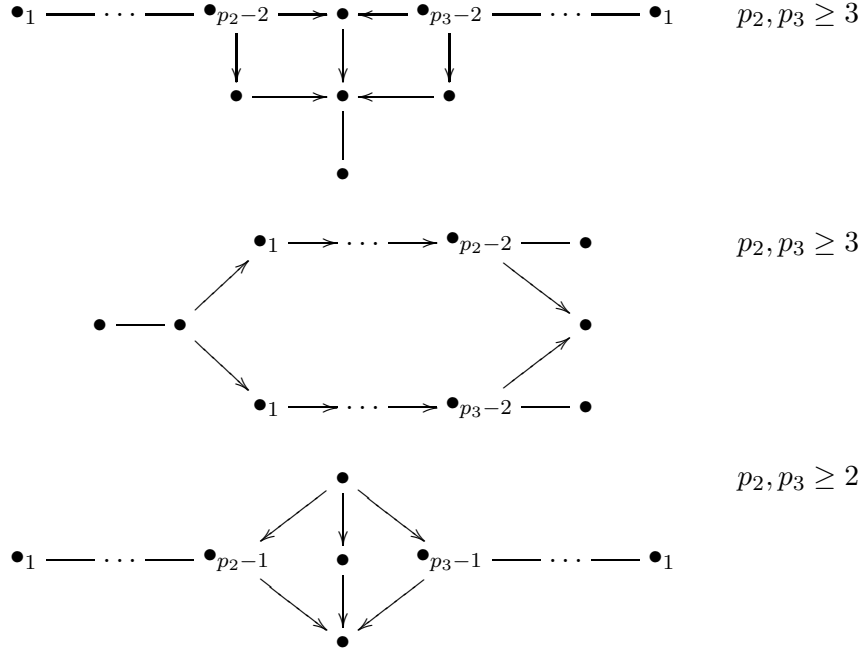
\[
\begin{array}{ccc}
\xymatrix@=1.5pc{
{\bullet_1} \ar@{-}[r] &
{\ldots} \ar@{-}[r] &
{\bullet_{p_2-2}} \ar[r] \ar[d] &
{\bullet} \ar[d] &
{\bullet_{p_3-2}} \ar[l] \ar[d] &
{\ldots} \ar@{-}[l] &
{\bullet_1} \ar@{-}[l] \\
& & {\bullet} \ar[r] &
{\bullet} \ar@{-}[d] &
{\bullet} \ar[l] \\
& & & {\bullet} \\
}
& &
p_2, p_3 \geq 3
\\
\\
\xymatrix@=1.5pc{
& & {\bullet_1} \ar[r] & {\ldots} \ar[r] &
{\bullet_{p_2-2}} \ar@{-}[r] \ar[dr] & {\bullet} \\
{\bullet} \ar@{-}[r] & {\bullet} \ar[ru] \ar[rd]
& & & & {\bullet} \\
& & {\bullet_1} \ar[r] & {\ldots} \ar[r] &
{\bullet_{p_3-2}} \ar@{-}[r] \ar[ur] & {\bullet}
}
& &
p_2, p_3 \geq 3
\\
\\
\xymatrix@=1.5pc{
& & & {\bullet} \ar[dl] \ar[d] \ar[dr] \\
{\bullet_1} \ar@{-}[r] & {\ldots} \ar@{-}[r] &
{\bullet_{p_2-1}} \ar[dr] & {\bullet} \ar[d] &
{\bullet_{p_3-1}} \ar[dl] \ar@{-}[r] &
{\ldots} \ar@{-}[r] & {\bullet_1} \\
& & & {\bullet}
}
& &
p_2, p_3 \geq 2
\end{array}
\]
where edges without arrows can be oriented arbitrarily.
\end{rem}

\subsection{Remarks on path algebras}
\label{ssec:quivers}
%%%%%%%%%%%%%%%%%%%%%%%%%%%%%%%%%%%%%
Let $Q$ be a finite quiver without oriented cycles. The set of its
vertices $Q_0$ has a natural partial order defined by $x \leq y$ for
two vertices $x, y \in Q_0$ if there exists an oriented path from $x$
to $y$. When $Q$ has the property that any two vertices are connected
by at most one oriented path, the path algebra $kQ$ is isomorphic to
the incidence algebra of the poset $(Q_0, \leq)$. A partial converse is
given by the following lemma.

\begin{lemma}[\protect{\cite[Theorem~4.2]{Mitchell68}}] \label{l:gldim2}
Let $X$ be a poset. Then $\gldim kX \leq 1$ if and only if any two
points in the Hasse diagram of $X$ are connected by at most one path.
\end{lemma}

For two quivers $Q$ and $Q'$, we denote $Q \sim Q'$ if $Q'$ can be
obtained from $Q$ by applying a sequence of BGP reflections (at sources
or sinks), see~\cite[(I.5.7)]{Happel88}. Since BGP reflections preserve
the derived equivalence class, we conclude that if $Q$ is a quiver such
that $Q \sim Q'$ for a quiver $Q'$ having the property that any two
vertices are connected by at most one oriented path, then the path
algebra $kQ$ is derived equivalent to an incidence algebra of a poset.

\subsection{The case $t=2$}
%%%%%%%%%%%%%%%%%%%%%%%%%%%%%%%
When the number of weights is at most 2, the corresponding canonical
algebra is a path algebra of a quiver, however there are two distinct
paths from the source $0$ to the sink $\omega$.

When the weight type is $(p_1, p_2)$ (with $p_1, p_2 \geq 2$), we can
overcome this problem by applying a BGP reflection at the sink
$\omega$. The resulting quiver is shown below,
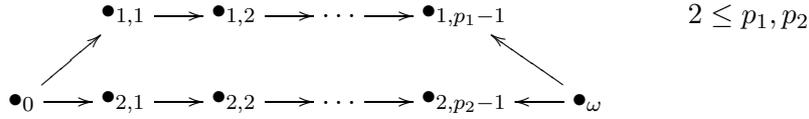
\begin{align*}
\xymatrix@=1.5pc{
& {\bullet_{1,1}} \ar[r] & {\bullet_{1,2}} \ar[r]
& {\ldots} \ar[r] & {\bullet_{1,p_1-1}} \\
{\bullet_0} \ar[ur] \ar[r]
& {\bullet_{2,1}} \ar[r] & {\bullet_{2,2}} \ar[r]
& {\ldots} \ar[r] & {\bullet_{2,p_2-1}}
& {\bullet_{\omega}} \ar[ul] \ar[l]
}
& & 2 \leq p_1, p_2
\end{align*}
and its path algebra is an incidence algebra derived equivalent to the
canonical algebra of type $(p_1,p_2)$.

\begin{figure}
\[
\xymatrix@=1.5pc{
{\bullet_0} \ar[rrrrr] \ar[dr] & & & & & {\bullet_p} \\
& {\bullet_1} \ar[r] & {\bullet_2} \ar[r] & {\ldots} \ar[r] &
{\bullet_{p-1}} \ar[ur] }
\]
\caption{The quiver $\widetilde{A}_{1,p}$ -- not derived equivalent to
any incidence algebra.} \label{f:A1p}
\end{figure}
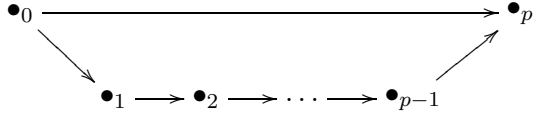

In the remaining case, where the weight type is either $(p)$ or $()$,
the corresponding canonical algebra equals the path algebra of the
quiver $\widetilde{A}_{1,p}$ drawn in Figure~\ref{f:A1p}, and we show
the following.

\begin{prop}
Let $p \geq 1$ and let $k$ be algebraically closed. Then there is no
poset whose incidence algebra is derived equivalent to the path algebra
of the quiver $\widetilde{A}_{1,p}$.
\end{prop}
\begin{proof}
Assume that there exists a poset $X$ such that $kX$ is derived
equivalent to the path algebra of $\widetilde{A}_{1,p}$.

If $\gldim kX \leq 1$, then by Lemma~\ref{l:gldim2}, the algebra $kX$
equals the path algebra of its Hasse diagram $Q$, thus by
\cite[(I.5.7)]{Happel88}, $Q$ is obtained from $\widetilde{A}_{1,p}$ by
a sequence of BGP reflections. But this is impossible since the only
possible reflections are at $0$ and $p$, and they give quivers
isomorphic to $\widetilde{A}_{1,p}$. However, $\widetilde{A}_{1,p}$ is
not the Hasse diagram of any poset.

Hence $\gldim kX \geq 2$. By Lemma~\ref{l:gldim2}, there exists at
least one commutativity relation in the quiver $Q$, so that $kX$ is not
a gentle algebra. This is again impossible since the path algebra of
$\widetilde{A}_{1,p}$ is gentle and the property of an algebra being
gentle is invariant under derived
equivalence~\cite{SchroerZimmermann03}. Alternatively, one can use the
explicit characterization in~\cite{AssemSkowronski87} of iterated
tilted algebras of type $\tilde{A}_n$.
\end{proof}

%%\bibliographystyle{acm}
%%\bibliography{canonical}

\def\cprime{$'$}

\end{document}